\documentclass[11pt]{amsart}
\usepackage[margin=1.1in]{geometry}
\usepackage{amsmath,amssymb,amsthm}
\usepackage{tikz}
\usepackage{booktabs}
\usepackage{xcolor}
\usepackage{enumitem}
\usepackage[colorlinks=true,linkcolor=blue,citecolor=blue,urlcolor=blue]{hyperref}

\newtheorem{theorem}{Theorem}[section]
\newtheorem{lemma}[theorem]{Lemma}
\newtheorem{proposition}[theorem]{Proposition}

\newtheorem{conjecture}[theorem]{Conjecture}
\theoremstyle{definition}
\newtheorem{definition}[theorem]{Definition}
\newtheorem{problem}[theorem]{Problem}
\theoremstyle{remark}
\newtheorem{remark}[theorem]{Remark}

\newcommand{\J}{J}
\newcommand{\GP}{\mathrm{GP}}

\title[Flower snarks and generalized Petersen graphs]{The cop number of the flower snarks, and the
generalized Petersen graphs with cop number two}
\author{Aarav Ghai}
\address{San Jos\'e State University}
\email{aaravmanav.ghai@sjsu.edu}
\date{\today}

\begin{document}

\begin{abstract}
We study two questions of Morris, Runte and Skelton (Australas.\ J.\ Combin.\ 83 (2022)) on the
game of Cops and Robbers.

First, we show that the flower snark $\J_n$ has cop number $3$ for every $n\ge 5$, and that $\J_3$
has cop number $2$. The upper bound comes from a three-cop strategy, which we call the
\emph{rolling wall}. It works on every connected graph whose vertices can be arranged in a cycle of
``columns'', each dominated by a single vertex. Because only one cop moves on each turn, the same
strategy also shows that the lazy cop number of $\J_n$ is $3$ for $n\ge 5$.

Second, we prove the conjecture of Morris, Runte and Skelton that the only generalized Petersen
graphs with cop number $2$ are $\GP(n,1)$ and the four graphs $\GP(6,2)$, $\GP(8,2)$, $\GP(9,3)$
and $\GP(12,3)$. The new ingredient is an explicit robber strategy showing that two cops do not
suffice on $\GP(3k,k)$ and $\GP(4k,k)$ for $k\ge 4$.

The theorem on flower snarks and all of the new lower bounds have been formally verified in the Lean~4 proof assistant.
\end{abstract}

\maketitle

\section{Introduction}\label{sec:intro}

In the game of \emph{Cops and Robbers}, introduced by Nowakowski and Winkler~\cite{NowakowskiWinkler1983}, a team of cops tries to catch a robber on a finite graph. Both sides move along
edges and see everything. The \emph{cop number} $c(G)$ is the least number of cops that can
guarantee a capture. Aigner and Fromme~\cite{AignerFromme1984} showed that graphs of girth at least
$5$ and minimum degree $\delta$ have $c(G)\ge\delta$. For background see the book of Bonato and
Nowakowski~\cite{BonatoNowakowski2011}.

The generalized Petersen graphs $\GP(n,k)$ are cubic, and they have been a testing ground for the
game. Ball, Bell, Guzman, Hanson-Colvin and Schonsheck~\cite{Ball2017} proved that
$c(\GP(n,k))\le 4$ for all $n,k$. Morris, Runte and Skelton~\cite{MorrisRunteSkelton2022} showed
that the bound $4$ is attained by most generalized Petersen graphs of girth $8$. They ended their
paper with several open questions. Two of them are answered here. The first concerns another
family of cubic graphs, the \emph{flower snarks} of Isaacs~\cite{Isaacs1975}.

\begin{problem}[{\cite[Problem~6.6]{MorrisRunteSkelton2022}}]\label{prob:flower}
What is the cop number of the flower snark $\J_n$?
\end{problem}

\begin{theorem}\label{thm:flower}
$c(\J_3)=2$, and $c(\J_n)=3$ for every $n\ge 5$. Moreover, the lazy cop number of $\J_n$ is $3$ for
every $n\ge 5$.
\end{theorem}

(Isaacs's flower snarks are the graphs $\J_n$ with $n$ odd, since only those are snarks. The
theorem holds for every $n\ge 5$.)

An earlier public, AI-generated exploration on ProbXiv~\cite{ProbXivFlower2026} also states
$c(\J_n)=3$ for odd $n\ge5$ and sketches a guarded-column strategy similar to the one below.
That entry reports no human examination. Our flower-snark treatment additionally proves the
lazy-cop statement and the even-$n$ extension, states the general cyclic-decomposition lemma,
and includes formal verification. We cite the entry as a prior public claim, not as a
peer-reviewed proof.

The second question is a conjecture stated in the same paper~\cite[p.~221]{MorrisRunteSkelton2022}.
A generalized Petersen graph is never cop-win, and it can have cop number $2$ only if its girth is
$3$ or $4$. All such graphs with $n\le 40$ are known~\cite{Ball2017}.

\begin{conjecture}[Morris, Runte, Skelton~\cite{MorrisRunteSkelton2022}]\label{conj:mrs}
The only generalized Petersen graphs with cop number $2$ are $\GP(n,1)$ for $n\ge3$ and
$\GP(3k,k)$, $\GP(4k,k)$ for $k\in\{2,3\}$.
\end{conjecture}

\begin{theorem}\label{thm:gp}
Conjecture~\ref{conj:mrs} is true. Equivalently, if $k\ge 4$ then
$c(\GP(3k,k))\ge 3$ and $c(\GP(4k,k))\ge 3$.
\end{theorem}

\subsection*{Methods}
The lower bound $c(\J_n)\ge 3$ follows from the Aigner--Fromme argument once one knows that $\J_n$
has girth at least $5$ for $n\ge5$. The upper bound is the substance of
Theorem~\ref{thm:flower}. In $\J_n$, the vertex $a_i$ is adjacent to $b_i,c_i,d_i$, so it dominates
the ``column'' $\{a_i,b_i,c_i,d_i\}$, and every edge joins two equal or consecutive columns. A cop
on $a_i$ is therefore a wall the robber cannot cross. Two such walls confine the robber, and a
third cop moves the walls closer together (Lemma~\ref{lem:wall}).

For Theorem~\ref{thm:gp} the girth-$3$ and girth-$4$ cases are $\GP(3k,k)$ and $\GP(4k,k)$. There
the Aigner--Fromme argument fails, and we replace it by an explicit set of safe positions for the
robber, which he can stay inside forever (Section~\ref{sec:gp}).

\subsection*{Verification}
\begin{itemize}
\item The statement ``$c(\J_n)=3$ for all $n\ge 5$'' has been formally verified in Lean~4 with
  Mathlib~\cite{Lean4,Mathlib2020}, for every $n$, without computer search
  (Section~\ref{sec:lean}).
\item So has the statement ``$c(\GP(n,k))\ge 3$ unless $k=1$ or $(n,k)$ is one of the four listed
  pairs''.
\item All the results were also checked by exhaustive computation for small parameters, using two
  independent solvers (Section~\ref{sec:computation}).
\end{itemize}

\section{Preliminaries}\label{sec:prelim}

All graphs are finite and simple. $N(x)$ is the neighbourhood of $x$ and $N[x]=N(x)\cup\{x\}$ is the
closed neighbourhood. We say that $x$ \emph{dominates} $y$ if $y\in N[x]$.

\subsection{The game}
A game with $k$ cops on a graph $G$ goes as follows.
\begin{enumerate}[label=(\arabic*)]
\item The cops choose $k$ vertices (not necessarily distinct).
\item The robber then chooses a vertex.
\item The cops and the robber then move alternately, starting with the cops. In a move, every cop
  (respectively the robber) stays put or moves to an adjacent vertex.
\end{enumerate}
The cops win if at some point a cop occupies the robber's vertex. The cop number $c(G)$ is the least
$k$ such that $k$ cops have a strategy that wins against every robber. In the \emph{lazy} variant at
most one cop moves on each cop turn~\cite{OffnerOjakian2014}. The \emph{lazy cop number} $c_L(G)$
is defined in the same way. Clearly $c(G)\le c_L(G)$.

A \emph{position} $(C,r)$ consists of the multiset $C$ of cop vertices and the robber's vertex $r$,
with the cops about to move. To show that $k$ cops do not suffice we use the following standard
principle.

\begin{lemma}\label{lem:safe}
Let $\mathcal S$ be a set of positions for $k$ cops such that:
\begin{enumerate}[label=(\alph*)]
\item \textup{(placement)} for every $C$ there is an $r$ with $(C,r)\in\mathcal S$;
\item \textup{(safety)} if $(C,r)\in\mathcal S$ then no cop of $C$ dominates $r$;
\item \textup{(closure)} if $(C,r)\in\mathcal S$ and $C'$ is obtained from $C$ by a cop move, then
  there is $r'\in N[r]$ with $(C',r')\in\mathcal S$.
\end{enumerate}
Then the robber can evade $k$ cops forever, so $c(G)>k$.
\end{lemma}

\begin{proof}
The robber starts at a vertex given by (a) and then always replies as in (c). By (b), when the cops
are about to move, none of them is on or next to the robber. So the cops' move cannot reach the
robber. After the robber's reply the position is again in $\mathcal S$, so by (b) the robber is not
standing on a cop.
\end{proof}

\begin{theorem}[Aigner--Fromme~\cite{AignerFromme1984}]\label{thm:af}
If $G$ has girth at least $5$ and minimum degree at least $3$, then $c(G)\ge 3$.
\end{theorem}

For completeness we recall the proof in the form we shall reuse.

\begin{proof}
Let $\mathcal S$ be the set of positions in which no cop dominates the robber.

\emph{Closure.} Suppose $(C,r)\in\mathcal S$ and the cops move to $C'$. No cop of $C'$ is on $r$.
Let $c\ne r$ be a vertex. If $N[c]$ contained two neighbours $x\ne y$ of $r$, there are two cases.
\begin{itemize}
\item If $c\in\{x,y\}$, then $r,x,y$ span a triangle.
\item Otherwise $r\,x\,c\,y$ is a $4$-cycle.
\end{itemize}
Both are impossible. So each cop dominates at most one of the at least three neighbours of $r$. The
robber moves to an undominated neighbour, and the new position is in $\mathcal S$.

\emph{Placement.} Here we need a vertex dominated by neither cop. For a cubic graph with more than
$8$ vertices one exists, since two closed neighbourhoods contain at most $8$ vertices. (For $\J_n$
see also Lemma~\ref{lem:flowerfacts}(c).)

Lemma~\ref{lem:safe} now applies with $k=2$.
\end{proof}

\section{Flower snarks}\label{sec:flower}

\begin{definition}[\cite{Isaacs1975}]\label{def:flower}
Let $n\ge 3$. The graph $\J_n$ has the $4n$ vertices $a_i,b_i,c_i,d_i$ with $0\le i\le n-1$, and the
following edges (indices of $b$ taken modulo $n$):
\begin{itemize}
  \item $a_ib_i$, $a_ic_i$, $a_id_i$ for every $i$ (the \emph{claws});
  \item $b_ib_{i+1}$ for every $i$ (the \emph{inner $n$-cycle});
  \item $c_ic_{i+1}$ and $d_id_{i+1}$ for $0\le i\le n-2$, together with $c_{n-1}d_0$ and
        $d_{n-1}c_0$ (the \emph{outer $2n$-cycle} $c_0c_1\cdots c_{n-1}d_0d_1\cdots d_{n-1}c_0$).
\end{itemize}
For odd $n$, $\J_n$ is the flower snark.
\end{definition}

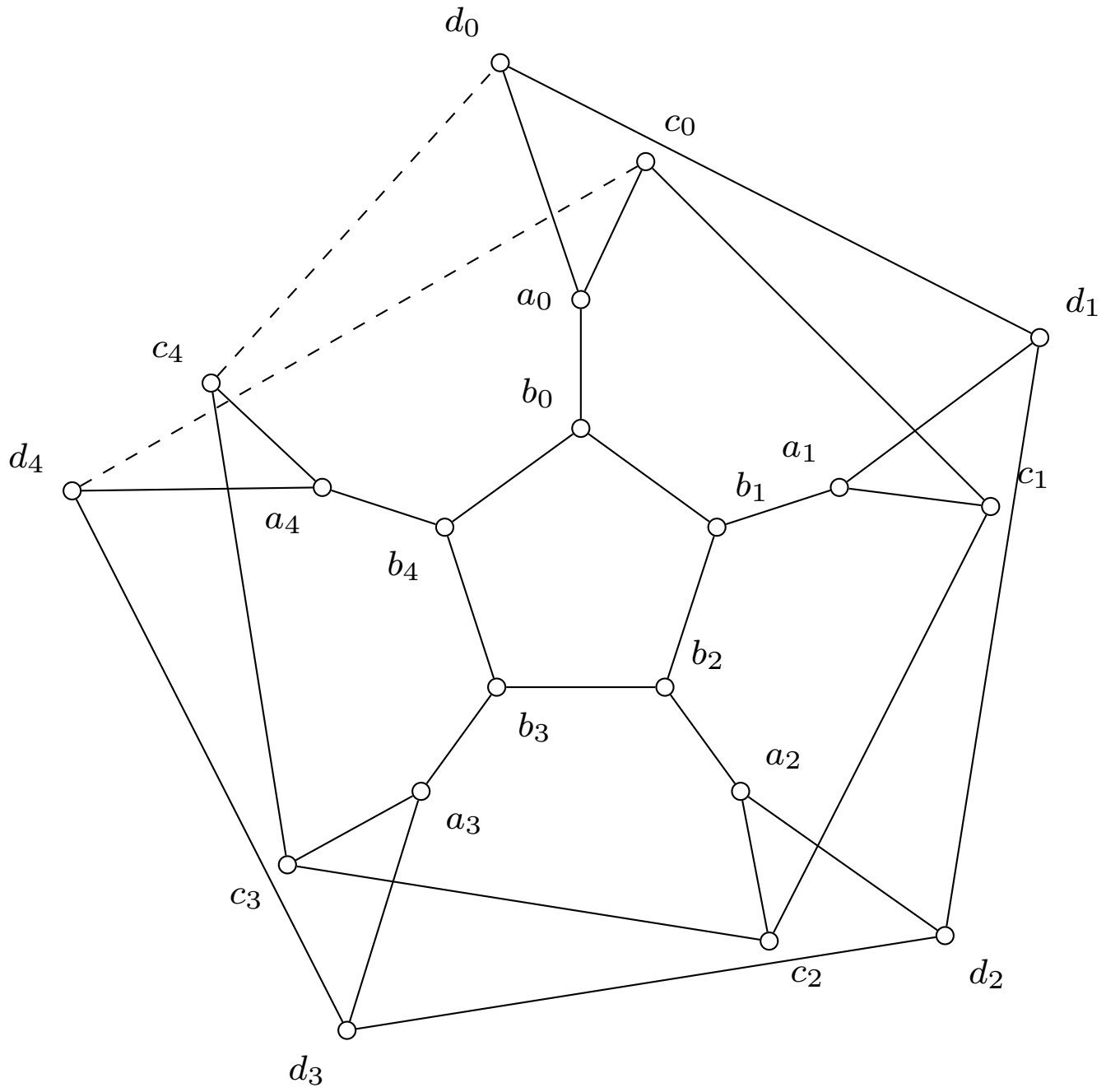
\begin{figure}[ht]
\centering
\begin{tikzpicture}[scale=1.15,
  v/.style={circle,draw,fill=white,inner sep=1.4pt},
  lbl/.style={font=\scriptsize}]
  \foreach \i in {0,...,4} {
    \pgfmathsetmacro{\ang}{90-72*\i}
    \node[v,label={[lbl]\ang+63:$b_{\i}$}] (b\i) at (\ang:1.0) {};
    \node[v,label={[lbl]\ang+90:$a_{\i}$}] (a\i) at (\ang:1.9) {};
    \node[v,label={[lbl]\ang-9:$c_{\i}$}] (c\i) at (\ang-9:2.9) {};
    \node[v,label={[lbl]\ang+9:$d_{\i}$}] (d\i) at (\ang+9:3.6) {};
  }
  \foreach \i in {0,...,4} {
    \pgfmathtruncatemacro{\j}{mod(\i+1,5)}
    \draw (a\i)--(b\i) (a\i)--(c\i) (a\i)--(d\i) (b\i)--(b\j);
  }
  \foreach \i in {0,...,3} {
    \pgfmathtruncatemacro{\j}{\i+1}
    \draw (c\i)--(c\j) (d\i)--(d\j);
  }
  \draw[dashed] (c4)--(d0) (d4)--(c0);
\end{tikzpicture}
\caption{The flower snark $\J_5$. The dashed edges $c_4d_0$ and $d_4c_0$ join the $c$-path and the
$d$-path into one $10$-cycle. Column $i$ is the claw $\{a_i,b_i,c_i,d_i\}$.}
\label{fig:J5}
\end{figure}

The \emph{column} of a vertex $x_i\in\{a_i,b_i,c_i,d_i\}$ is $\pi(x_i)=i\in\mathbb Z_n$.

\begin{lemma}\label{lem:flowerfacts}
Let $n\ge 3$.
\begin{enumerate}[label=(\alph*)]
\item $N[a_i]$ is exactly column $i$, and every edge of $\J_n$ joins two vertices whose columns are
  equal or differ by $\pm1$ modulo $n$.
\item If $n\ge 5$, then $\J_n$ has girth at least $5$ (in fact $5$ for $n=5$ and $6$ for $n\ge 7$).
\item Every closed neighbourhood contains at most one vertex of $\{a_0,\dots,a_{n-1}\}$.
\end{enumerate}
\end{lemma}

\begin{proof}
(a) The first claim is immediate. Claws stay inside a column. The edges $b_ib_{i+1}$, $c_ic_{i+1}$
and $d_id_{i+1}$ join consecutive columns, and the two twist edges join column $n-1$ to column $0$.

(c) The $a$-vertices are pairwise non-adjacent, and each other vertex has exactly one
$a$-neighbour.

(b) We show that no cycle of length at most $4$ passes through an $a$-vertex. Cycles avoiding the
$a$-vertices lie in the inner $n$-cycle or the outer $2n$-cycle, so they are long.

Suppose a short cycle passes through $a_i$, with cycle-neighbours $x\ne y$ in $\{b_i,c_i,d_i\}$.
\begin{itemize}
\item These three vertices are pairwise non-adjacent, so the cycle is not a triangle.
\item For a $4$-cycle, $x$ and $y$ would need a common neighbour $z\ne a_i$. Apart from $a_i$, the
  neighbours of $b_i$ are $b$-vertices, while those of $c_i,d_i$ lie on the outer cycle. So
  $\{x,y\}=\{c_i,d_i\}$, and $z$ is on the outer cycle.
\item On the outer $2n$-cycle, $c_i$ and $d_i$ are at positions $i$ and $n+i$. A common neighbour
  would require $n\equiv\pm2\pmod{2n}$, which is impossible for $n\ge 3$.
\item Finally, the inner and outer cycles have lengths $n\ge5$ and $2n\ge 10$.
\end{itemize}
The exact values of the girth were confirmed by computer for $n\le 41$.
\end{proof}

\subsection{The lower bound}

\begin{proposition}\label{prop:lower}
For $n\ge 5$, two cops cannot catch the robber on $\J_n$. Hence $c(\J_n)\ge 3$.
\end{proposition}

\begin{proof}
$\J_n$ is cubic with girth at least $5$ by Lemma~\ref{lem:flowerfacts}(b), so
Theorem~\ref{thm:af} applies. For placement, two cops dominate at most two of $a_0,a_1,a_2$ by
Lemma~\ref{lem:flowerfacts}(c), and the robber starts on the third.
\end{proof}

The hypothesis $n\ge5$ is needed. $\J_3$ contains the triangle $b_0b_1b_2$, and $\J_4$ contains
the $4$-cycle $b_0b_1b_2b_3$.

\subsection{The upper bound: the rolling wall}

We prove the upper bound for a class of graphs that contains every $\J_n$.

\begin{definition}\label{def:decomp}
A \emph{cyclic dominated decomposition} of a graph $G$ into $n$ columns is a map
$\pi\colon V(G)\to\{0,\dots,n-1\}$ together with vertices $a_0,\dots,a_{n-1}$ such that:
\begin{enumerate}[label=(\roman*)]
\item each vertex $v$ is dominated by $a_{\pi(v)}$;
\item every edge $uv$ satisfies $|\pi(u)-\pi(v)|\le 1$, or $\{\pi(u),\pi(v)\}=\{0,n-1\}$.
\end{enumerate}
\end{definition}

By Lemma~\ref{lem:flowerfacts}(a), the columns of $\J_n$ form such a decomposition. Its key
consequence is the following.

\emph{A cop on $a_j$ is a wall.} If the robber ever enters column $j$, he is next to that cop,
who captures him on the following move. Hence a robber who wants to survive never enters column $j$.

\begin{lemma}[Rolling wall]\label{lem:wall}
If $G$ is connected and has a cyclic dominated decomposition, then three cops win on $G$. They win
even in the lazy game, so $c(G)\le c_L(G)\le 3$.
\end{lemma}

\begin{proof}
Call the cops $W$, $F$ and $M$ (back wall, front wall, mover). All three start on $a_0$.
\begin{itemize}
\item Throughout, if the robber is dominated by some cop at the start of a cop turn, that cop moves
  onto him and the game ends.
\item Otherwise exactly one cop moves, as described below. So the strategy is lazy.
\end{itemize}
We keep an integer $f$, initially $0$, and maintain the following.

\smallskip
\noindent\textbf{Invariant.} At each cop turn, $W$ is on $a_0$, $F$ is on $a_f$, and the robber is
in a column in $\{f+1,\dots,n-1\}$.
\smallskip

\emph{The invariant holds at the start.} After the robber places, either he is in column $0$ and is
caught at once, or he is in a column in $\{1,\dots,n-1\}$.

\emph{The robber cannot leave the interval.} Suppose the invariant holds and only $M$ moves. If the
robber steps out of $\{f+1,\dots,n-1\}$, then by (ii) he enters column $f$ or column $0$. There he
is dominated by $F$ or $W$, by (i), and is captured on the next turn.

\emph{A phase.} By connectivity, $M$ walks along a path to $a_{f+1}$, one edge per turn, while $W$
and $F$ stand still. During the walk the invariant is preserved, by the previous paragraph. When
$M$ arrives on $a_{f+1}$ there are two cases.
\begin{itemize}
\item If the robber is in column $f+1$, $M$ dominates him and captures him.
\item Otherwise he is in a column in $\{f+2,\dots,n-1\}$.
\end{itemize}
We now rename the cops: $M$ becomes the new front wall, $F$ becomes the new mover, and $f$
increases by one. The invariant holds again.

\emph{Termination.} Each phase is finite and increases $f$. When $f=n-1$ the interval
$\{f+1,\dots,n-1\}$ is empty, so the robber must already have been caught.
\end{proof}

\begin{proof}[Proof of Theorem~\ref{thm:flower}]
For $n\ge 5$, Proposition~\ref{prop:lower} gives $c(\J_n)\ge3$. Lemmas~\ref{lem:wall}
and~\ref{lem:flowerfacts}(a) give $c_L(\J_n)\le 3$. Since $c\le c_L$, both are equal to $3$.

For $n=3$, exhaustive computation (Section~\ref{sec:computation}) shows that one cop loses and two
cops win on $\J_3$. Two independent programs agree on this.
\end{proof}

\begin{remark}\label{rem:wall}
The strategy is wasteful. Each phase moves $M$ from $a_{f-1}$ to $a_{f+1}$ in four steps (via
$b_{f-1},b_f,b_{f+1}$), and its worst-case capture time is exactly $4n-5$ (computed for odd
$n\le 151$). The optimal capture time with three cops is much smaller, roughly $n/2$
(Table~\ref{tab:data}).

The decomposition idea is close to standard ``guard a separator'' arguments; compare the proof
of $c(G)\le \mathrm{tw}(G)/2+1$ in~\cite{JoretKaminskiTheis2010}. We have not found
Lemma~\ref{lem:wall} stated in the literature, but we do not claim the technique is new.
\end{remark}

\section{Generalized Petersen graphs with cop number two}\label{sec:gp}

For $n\ge 3$ and $1\le k<n/2$, the generalized Petersen graph $\GP(n,k)$~\cite{Watkins1969} has
outer vertices $u_j$ and inner vertices $v_j$ ($j\in\mathbb Z_n$). Its edges are $u_ju_{j+1}$,
$u_jv_j$ and $v_jv_{j+k}$. Thus
\[
N[u_j]=\{u_{j-1},u_j,u_{j+1},v_j\},\qquad N[v_j]=\{v_j,u_j,v_{j+k},v_{j-k}\}.
\]

\begin{lemma}[cf.~\cite{Ball2017}]\label{lem:gpgirth}
Let $2\le k<n/2$. Then $\GP(n,k)$ has girth at most $4$ if and only if $n=3k$ or $n=4k$.
\end{lemma}

\begin{proof}
Every cycle uses an even number of spokes $u_jv_j$. Consider the short cycles by the number of
spokes they use.
\begin{itemize}
\item \emph{No spokes.} The cycle lies in the outer $n$-cycle, so its length is $n\ge5$; or it lies
  in the inner graph. A short inner cycle needs $3k\equiv 0$ or $4k\equiv 0\pmod n$. Since $2k<n$,
  this means $n=3k$ or $n=4k$.
\item \emph{Two spokes, length $4$.} The cycle is $u_ju_{j\pm1}v_{j\pm1}v_j$, which needs
  $k\equiv\pm1\pmod n$. That is excluded.
\item \emph{Two spokes, length $3$.} This is impossible.
\end{itemize}
\end{proof}

\begin{theorem}\label{thm:gprobber}
If $k\ge4$, two cops cannot catch the robber on $\GP(3k,k)$ or on $\GP(4k,k)$.
\end{theorem}

We first show how Theorem~\ref{thm:gp} follows.

\begin{proof}[Proof of Theorem~\ref{thm:gp} from Theorem~\ref{thm:gprobber}]
The listed graphs have cop number $2$.
\begin{itemize}
\item For $\GP(n,1)$, $\GP(6,2)$ and $\GP(8,2)$ this is in~\cite{Ball2017} (by hand).
\item For $\GP(9,3)$ and $\GP(12,3)$ it is also in~\cite{Ball2017} (by computer). We
  confirmed all four of the small cases with two independent solvers.
\end{itemize}
Now let $2\le k<n/2$ with $(n,k)$ not listed. There are two cases.
\begin{itemize}
\item If $n\notin\{3k,4k\}$, then $\GP(n,k)$ is cubic of girth at least $5$ by
  Lemma~\ref{lem:gpgirth}, so $c\ge3$ by Theorem~\ref{thm:af}. Placement holds since $2n>8$.
\item Otherwise $k\ge 4$, and $c\ge3$ by Theorem~\ref{thm:gprobber}.
\end{itemize}
\end{proof}

\subsection{The safe set}
Fix $n=3k$ or $n=4k$ with $k\ge4$. For an index $i$ put $U_i=\{u_{i-1},u_i,u_{i+1}\}$ and
\[
B_i=\{u_{i-2},u_{i-1},u_{i+1},u_{i+2},v_{i-1},v_{i+1}\}.
\]
A cop that can step into $U_i$ in one move and is not in $N[v_i]$ lies in $B_i$. We also define
\[
A_i=\begin{cases}\{u_{i-k},u_{i+k}\} & \text{if } n=3k,\\
\{v_{i+2k},u_{i+2k}\} & \text{if } n=4k.\end{cases}
\]
Let $T_i$ denote the inner triangle $\{v_i,v_{i+k},v_{i-k}\}$ when $n=3k$, and let
$w_i=v_{i+2k}$ be the antipode of $v_i$ on its inner $4$-cycle when $n=4k$. Then $A_i$ is the set
of vertices outside $N[v_i]$ from which a cop can step onto a vertex dominating both inner
neighbours of $v_i$: a vertex of $T_i$, or $w_i$.

Let $\mathcal S$ be the set of positions $(\{c_1,c_2\},r)$ such that:
\begin{enumerate}[label=(S\arabic*)]
\item no cop dominates $r$; and
\item if $r=v_i$ is an inner vertex, it is not the case that one cop is in $A_i$ and the other is
  in $B_i$.
\end{enumerate}
In words, a robber on $v_i$ is threatened only if one cop guards the exit to the outer rim while the
other guards both inner exits. (S2) forbids exactly this combination.

\subsection{Index facts}
Put $D_i=T_i$ if $n=3k$, and $D_i=\{v_i,w_i\}$ if $n=4k$. All indices are modulo $n$, and $e$ ranges
over $\{-1,0,1\}$. Each fact below says that some small integer combination of $k$ is not
$\equiv0\pmod n$. Each is checked directly and holds under the stated hypothesis.
\begin{enumerate}[label=(F\arabic*)]
\item $v_{i\pm1}\notin N[v_i]$, and $u_{i+e}\notin N[v_{i\pm k}]$ \hfill ($k\ge2$).
\item $A_i\cap B_i=\emptyset$ \hfill ($k\ge 3$).
\item $N[v_{i+k}]\cap N[v_{i-k}]=D_i$, and $N[D_i]\setminus N[v_i]=A_i$ \hfill ($k\ge2$).
\item $u_{i+e}\notin B_{i+k}\cup B_{i-k}$ \hfill ($k\ge 4$).
\item $N[v_{i\pm k}]\setminus D_i$ is disjoint from $B_{i\mp k}$. This set is $\{u_{i\pm k}\}$ if
  $n=3k$, and $\{v_{i\pm k},u_{i\pm k}\}$ if $n=4k$. \hfill ($k\ge 3$).
\item $B_{i+k}\cap B_{i-k}=\emptyset$ \hfill ($k\ge5$ if $n=3k$; $k\ge3$ if $n=4k$).
\end{enumerate}
For example:
\begin{itemize}
\item (F4) fails only if $u_{i+e}=u_{i\pm k+\delta}$ with $|\delta|\le2$. That forces
  $k\equiv m\pmod n$ with $|m|\le 3$, which is impossible for $k\ge4$.
\item (F6) fails only if $2k\equiv m\pmod n$ with $|m|\le4$. For $n=3k$ this forces $k\le4$, and
  indeed (F6) fails for $\GP(12,4)$, where $u_{i+6}=u_{i-6}$.
\end{itemize}

\subsection{Proof of Theorem~\ref{thm:gprobber}}
We verify the three conditions of Lemma~\ref{lem:safe}. Safety is (S1).

\emph{Placement.}
\begin{itemize}
\item A cop on $u_m$ dominates three outer vertices, and a cop on $v_m$ dominates one.
\item Since $n>6$, some outer vertex is undominated. (S2) is vacuous there.
\end{itemize}

\emph{Closure.} Let $(C,r)\in\mathcal S$, and let the cops move to $C'$. We write $X'$ for the new
position of a cop $X$.

\smallskip\noindent\emph{Case 1: $r=u_i$.} No cop of $C'$ is on $u_i$.

If some vertex of $U_i$ is undominated by $C'$, the robber moves there.

Otherwise, since $n\ge5$, the only vertices other than $u_i$ that dominate two vertices of $U_i$ are
$u_{i-1}$ and $u_{i+1}$. So one cop is on, say, $u_{i+1}$, and the other dominates $u_{i-1}$
without being on $u_i$. Hence both cops lie in $B_i$. The robber moves to $v_i$.
\begin{itemize}
\item (S1) holds by (F1).
\item (S2) holds because $A_i\cap B_i=\emptyset$ (F2), so neither cop is in $A_i$.
\end{itemize}

\smallskip\noindent\emph{Case 2: $r=v_i$.} By (S1) no cop of $C$ is in $N[v_i]$, so no cop of
$C'$ is on $v_i$.

\emph{Case 2a.} If no cop of $C'$ is in $U_i$, then $u_i$ is undominated, and the robber moves to
$u_i$.

\emph{Case 2b.} Otherwise, let $X$ be a cop with $X'\in U_i$, and let $Y$ be the other cop.
\begin{itemize}
\item $X$ was in $N[U_i]\setminus N[v_i]=B_i$.
\item So, by (S2) for $(C,v_i)$, $Y$ was not in $A_i$. By (S1), $Y$ was not in $N[v_i]$ either.
\item By (F1) and (F4), $X'$ dominates neither $v_{i+k}$ nor $v_{i-k}$, and
  $X'\notin B_{i+k}\cup B_{i-k}$.
\end{itemize}

We claim that there is $t\in\{v_{i+k},v_{i-k}\}$ that $Y'$ does not dominate and with
$Y'\notin B_t$.
\begin{itemize}
\item First, $Y'\notin D_i$. Otherwise $Y\in N[D_i]\setminus N[v_i]=A_i$ by (F3).
\item Hence, by (F3), $Y'$ dominates at most one of $v_{i\pm k}$.
\item If $Y'$ dominates $v_{i\pm k}$, then $Y'\in N[v_{i\pm k}]\setminus D_i$. Take $t=v_{i\mp k}$,
  so $Y'\notin B_t$ by (F5).
\item If $Y'$ dominates neither, then by (F6) $Y'$ lies outside at least one of $B_{i+k}$ and
  $B_{i-k}$. Take $t$ to be the corresponding vertex.
\end{itemize}

The robber moves to $t$, which is legal since $t\in N[v_i]$.
\begin{itemize}
\item (S1) holds at $t$, since neither cop dominates it.
\item (S2) holds at $t$, since neither cop is in $B_t$.
\end{itemize}
This proves closure for $\GP(3k,k)$ with $k\ge5$ and for $\GP(4k,k)$ with $k\ge4$.

\emph{The case $\GP(12,4)$.} Here (F6) fails, and so does the claim in Case~2b. We use instead the safe set $\mathcal S_{12}$ obtained
from $\mathcal S$ by also forbidding, for a robber on $v_i$, the eight cop pairs
\[
\{u_{i\pm1},u_{i\pm5}\}\ (\text{four pairs}),\qquad
\{u_{i\pm1},u_{i+6}\},\qquad \{u_{i\pm1},v_{i+6}\}.
\]
A direct enumeration of all $300$ cop placements and $24$ robber vertices verifies conditions
(a)--(c) of Lemma~\ref{lem:safe} for $\mathcal S_{12}$. This finishes the proof. \qed

\begin{remark}
For every $k\le 40$, the set $\mathcal S$ (respectively $\mathcal S_{12}$) coincides exactly with
the set of positions from which the robber can evade two cops, as computed by an exhaustive solver.
The robber's replies described above were also replayed mechanically against every cop move for
$k\le 30$.
\end{remark}

\section{Computations}\label{sec:computation}

We computed the outcome of the game by retrograde analysis of the game graph, which gives the
exact winning region and the exact capture time. Two independent programs were used. They share no
code, and each builds the graphs from the definitions. They agree on every case that both
computed.

\begin{table}[ht]
\centering
\small
\begin{tabular}{@{}l*{10}{c}@{}}
\toprule
$n$ & 5 & 7 & 9 & 11 & 13 & 15 & 17 & 19 & 21 & 23 \\
capture time & 2 & 3 & 4 & 5 & 6 & 7 & 8 & 8 & 9 & 10 \\
\midrule
$n$ & 25 & 27 & 29 & 31 & 33 & 35 & 37 & 39 & 41 & \\
capture time & 11 & 11 & 13 & 13 & 14 & 15 & 16 & 16 & 18 & \\
\bottomrule
\end{tabular}
\caption{Optimal three-cop capture time on $\J_n$, measured in cop moves. For every odd $n\le 41$,
two cops lose and three cops win; $\J_3$ has cop number $2$.}
\label{tab:data}
\end{table}

\begin{itemize}
\item \emph{Flower snarks.} Two cops lose and three cops win on $\J_n$ for every odd $5\le n\le
  41$. One cop loses and two cops win on $\J_3$. The lazy cop number of $\J_n$ is $3$ for
  $n\in\{3,5,7,9\}$; note that $c_L(\J_3)=3>c(\J_3)$.
\item \emph{Generalized Petersen graphs.} Two cops lose on $\GP(3k,k)$ and $\GP(4k,k)$ for
  $4\le k\le 40$ (up to $320$ vertices), and three cops win for $k\le 20$. Two cops win for
  $k\in\{2,3\}$.
\end{itemize}

\section{Formal verification}\label{sec:lean}

The statement
\[
\forall n\ge 5,\quad \operatorname{copNumber}(\J_n)=3
\]
is proved in Lean~4 using Mathlib~\cite{Lean4,Mathlib2020}. The proof assistant reports that it
depends only on the standard axioms (propositional extensionality, choice and quotient soundness).
It uses no computer search.

\begin{itemize}
\item \emph{The game.} A position is cop-winning if it lies in the least set of positions that is
  closed under ``the cops can move so that every robber reply leads into the set''. Defined this
  way, a winning position means capture is forced in finitely many moves.
\item \emph{Robber wins.} These are certified by a safe set as in Lemma~\ref{lem:safe}. We prove
  that the two notions exclude each other.
\item \emph{The lower bound.} This is the argument of Theorem~\ref{thm:af}. The local facts for
  $\J_n$ are discharged by linear arithmetic, uniformly in $n$.
\item \emph{The upper bound.} This is Lemma~\ref{lem:wall}, proved for an abstract cyclic dominated
  decomposition. The proof is by induction on the number of remaining columns, with an inner
  induction on the walk taken by the mover.
\item \emph{Fidelity of the definition.} The formal definition of $\J_n$ follows
  Definition~\ref{def:flower} literally. Its edge set was exported and checked to be isomorphic to
  the graphs used in Section~\ref{sec:computation}.
\end{itemize}

The results of Section~\ref{sec:gp} are formalised in the same development. The formal statement
says that if $k\ge2$, $2k<n$ and $(n,k)\notin\{(6,2),(8,2),(9,3),(12,3)\}$, then every winning team
of cops on $\GP(n,k)$ has at least three cops. Hence $\operatorname{copNumber}(\GP(n,k))\ge3$.
\begin{itemize}
\item \emph{Girth at least $5$.} When $n\notin\{3k,4k\}$, the hypotheses of the Aigner--Fromme
  argument are again proved by linear arithmetic. The four graphs with $n<9$ are handled by direct
  evaluation.
\item \emph{$\GP(3k,k)$ and $\GP(4k,k)$.} The safe set $\mathcal S$ and the case analysis of
  Section~\ref{sec:gp} are formalised uniformly in $k$.
\item \emph{$\GP(12,4)$.} The certificate $\mathcal S_{12}$ is checked by the Lean kernel itself.
  This is a kernel evaluation, not compiled code, so it adds no axioms.
\item \emph{Axioms and fidelity.} All of these theorems depend only on the three standard axioms.
  The formal $\GP(n,k)$ was checked to have exactly the textbook edge set for several $(n,k)$.
\end{itemize}
The matching upper bounds ($c=2$ for the listed graphs, from~\cite{Ball2017}) are not formalised.

\section{Concluding remarks}\label{sec:conclusion}

\begin{enumerate}
\item \emph{Capture time.} The optimal capture time of three cops on $\J_n$ grows roughly like
  $n/2$, but irregularly (Table~\ref{tab:data}). What is it exactly?
\item \emph{Lazy cops.} Theorem~\ref{thm:flower} determines the lazy cop number of $\J_n$, giving a
  partial answer to the snark analogue of~\cite[Problem~6.3]{MorrisRunteSkelton2022}. The lazy cop
  number of generalized Petersen graphs remains open.
\item \emph{Dominated decompositions.} Lemma~\ref{lem:wall} applies to any graph that maps onto a
  cycle with dominated fibres. Which other snark families, or other cubic graphs, admit such a
  decomposition, possibly after contracting a bounded set?
\end{enumerate}

\subsection*{Acknowledgements}
Anthropic Claude assisted with proof-strategy exploration, manuscript drafting,
computational checks, and Lean formalization; the author takes responsibility for
the results.

\bibliographystyle{amsplain}
\bibliography{refs}

@article{MorrisRunteSkelton2022,
  author  = {Morris, Joy and Runte, Tigana and Skelton, Adrian},
  title   = {Most generalized {P}etersen graphs of girth 8 have cop number 4},
  journal = {Australasian Journal of Combinatorics},
  volume  = {83},
  number  = {2},
  pages   = {204--224},
  year    = {2022},
  url     = {https://ajc.maths.uq.edu.au/pdf/83/ajc_v83_p204.pdf}
}

@article{AignerFromme1984,
  author  = {Aigner, M. and Fromme, M.},
  title   = {A game of cops and robbers},
  journal = {Discrete Applied Mathematics},
  volume  = {8},
  number  = {1},
  pages   = {1--12},
  year    = {1984},
  doi     = {10.1016/0166-218X(84)90073-8}
}

@article{Ball2017,
  author  = {Ball, Taylor and Bell, Robert W. and Guzman, Jonathan and Hanson-Colvin, Madeleine and Schonsheck, Nikolas},
  title   = {On the cop number of generalized {P}etersen graphs},
  journal = {Discrete Mathematics},
  volume  = {340},
  number  = {6},
  pages   = {1381--1388},
  year    = {2017},
  doi     = {10.1016/j.disc.2016.10.004}
}

@article{Isaacs1975,
  author  = {Isaacs, Rufus},
  title   = {Infinite families of nontrivial trivalent graphs which are not {T}ait colorable},
  journal = {The American Mathematical Monthly},
  volume  = {82},
  number  = {3},
  pages   = {221--239},
  year    = {1975},
  doi     = {10.1080/00029890.1975.11993805}
}

@article{NowakowskiWinkler1983,
  author  = {Nowakowski, Richard and Winkler, Peter},
  title   = {Vertex-to-vertex pursuit in a graph},
  journal = {Discrete Mathematics},
  volume  = {43},
  number  = {2--3},
  pages   = {235--239},
  year    = {1983},
  doi     = {10.1016/0012-365X(83)90160-7}
}

@book{BonatoNowakowski2011,
  author    = {Bonato, Anthony and Nowakowski, Richard J.},
  title     = {The Game of Cops and Robbers on Graphs},
  series    = {Student Mathematical Library},
  volume    = {61},
  publisher = {American Mathematical Society},
  address   = {Providence, RI},
  year      = {2011},
  doi       = {10.1090/stml/061}
}

@article{OffnerOjakian2014,
  author  = {Offner, David and Ojakian, Kerry},
  title   = {Variations of cops and robber on the hypercube},
  journal = {Australasian Journal of Combinatorics},
  volume  = {59},
  number  = {2},
  pages   = {229--250},
  year    = {2014}
}

@article{JoretKaminskiTheis2010,
  author  = {Joret, Gwena{\"e}l and Kami{\'n}ski, Marcin and Theis, Dirk Oliver},
  title   = {The Cops and Robber game on graphs with forbidden (induced) subgraphs},
  journal = {Contributions to Discrete Mathematics},
  volume  = {5},
  number  = {2},
  year    = {2010},
  doi     = {10.55016/ojs/cdm.v5i2.62032}
}

@article{Watkins1969,
  author  = {Watkins, Mark E.},
  title   = {A theorem on {T}ait colorings with an application to the generalized {P}etersen graphs},
  journal = {Journal of Combinatorial Theory},
  volume  = {6},
  pages   = {152--164},
  year    = {1969},
  doi     = {10.1016/S0021-9800(69)80116-X}
}

@inproceedings{Lean4,
  author    = {de Moura, Leonardo and Ullrich, Sebastian},
  title     = {The {L}ean 4 Theorem Prover and Programming Language},
  booktitle = {Automated Deduction -- CADE 28},
  series    = {Lecture Notes in Computer Science},
  volume    = {12699},
  pages     = {625--635},
  publisher = {Springer},
  year      = {2021},
  doi       = {10.1007/978-3-030-79876-5_37}
}

@inproceedings{Mathlib2020,
  author    = {{The mathlib Community}},
  title     = {The {L}ean Mathematical Library},
  booktitle = {Proceedings of the 9th ACM SIGPLAN International Conference on Certified Programs and Proofs (CPP 2020)},
  pages     = {367--381},
  year      = {2020},
  doi       = {10.1145/3372885.3373824}
}

@misc{ProbXivFlower2026,
  key          = {ProbXiv},
  title        = {{AI}-generated exploration of the cop number of the flower snarks},
  howpublished = {ProbXiv, \url{https://probxiv.com/problems/most-generalized-petersen-graphs-of-girth-8-have-cop-number-4-3}},
  year         = {2026},
  note         = {Public comment dated 17 August 2026; the page states that no person has examined it}
}

\end{document}